\documentclass[11pt]{article}
\usepackage{indentfirst, latexsym, bm}
\usepackage{graphics}
\usepackage{amsmath}
\usepackage{amssymb}
\usepackage{amsfonts}
\usepackage{bbding}
\usepackage{dsfont}

\numberwithin{equation}{section} \topmargin=-1.5cm
\numberwithin{equation}{section}

\begin{document}

\newtheorem{theorem}{Theorem}[section]
\newtheorem{lemma}[theorem]{Lemma}
\newtheorem{proof}{Proof}

\title{EXTEND MEAN CURVATURE FLOW WITH FINITE INTEGRAL CURVATURE
\footnote{2000 Mathematics Subject Classification. 53C44; 53C21.
\newline\indent Research supported by the NSFC, Grant No. 10771187; the Trans-Century Training Programme Foundation for
Talents by the Ministry of Education of China; and the Natural
Science Foundation of Zhejiang Province, Grant No. 101037.
\newline \indent Keywords: Mean curvature flow, maximal existence
time, second fundamental form, integral curvature.}}
\author{H{\footnotesize ONG-WEI} X{\footnotesize U},\ \ F{\footnotesize EI} Y{\footnotesize E}\ {\footnotesize AND}\ E{\footnotesize N-TAO} Z{\footnotesize HAO} \\
}
\date{}
\maketitle


\begin{abstract}
In this note, we first prove that the solution of mean curvature
flow on a finite time interval $[0,T)$ can be extended over time $T$
if the space-time integration of the norm of the second fundamental
form is finite. Secondly, we prove that the solution of certain mean
curvature flow on a finite time interval $[0,T)$ can be extended
over time $T$ if the space-time integration of the mean curvature is
finite. Moreover, we show that these conditions are optimal in some
sense.

\end{abstract}



\section{Introduction}
Let $M$ be a complete $n$-dimensional manifold without boundary, and
let $F_t :M^n\rightarrow \mathbb{R}^{n+1}$ be a one-parameter family
of smooth hypersurfaces immersed in Euclidean space. We say that
$M_t = F_t(M)$ is a solution of the mean curvature flow if $F_t$
satisfies
\[
\left\{
\begin{array}{ccc}
\frac{\partial}{\partial t}F(x,t)&=&-H(x,t)\overrightarrow{\nu}(x,t)\\
F(x,0)&=&F_0(x),
\end{array}\right.\]
where $F(x,t)=F_t(x)$, $H(x,t)$ is the mean curvature,
$\overrightarrow{\nu}(x,t)$ is the unit outward normal vector, and
$F_{0}$ is some given initial hypersurface.

K. Brakke [1] studied the mean curvature flow from the view point of
geometric measure theory firstly. For the classical solution of the
mean curvature flow, G. Huisken (see [5], [6]) showed that for a
smooth complete initial hypersurface with bounded second fundamental
form the solution exists on a maximal time interval $[0, T)$,
$0<T\leq\infty$. If the closed initial hypersurface is convex, he
showed that in [6] the mean curvature flow will converge to a round
point in finite time. He also proved that if the second fundamental
form is uniformly bounded, then the mean curvature flow can be
extended.

By a blow up argument, N. \v{S}e\v{s}um [9] proved that if the Ricci
curvature is uniformly bounded on $M \times[0, T)$, then the Ricci
flow can be extended over $T$. In [10], B. Wang obtained some
integral conditions to extend the Ricci flow. A natural question is
that, what is the optimal condition for the mean curvature flow to
be extended? By a different method, we investigate the integral
conditions to extend the mean curvature flow. We will prove that the
mean curvature flow can be extended if
the integration of the norm of the second fundamental form is bounded. More preciously, we obtain the following\\\\
\textbf{Theorem 1.1.} \emph{Let $F_t:\ M^n\longrightarrow
\mathbb{R}^{n+1}$ be a solution of the mean curvature flow of closed
hypersurfaces on a finite time interval $[0,T)$. If
$$|| A||_{\alpha,M\times [0,T)}=
\left(\int^T_0\int_M | A |^\alpha d\mu
dt\right)^{\frac{1}{\alpha}}<+\infty,$$ for some $\alpha\geq n+2$, then this flow can be extended over time $T$.}\\

When the space-time integration of the mean curvature is finite and
the second fundamental tensor is bounded from below, we also prove
the following theorem.\\\\
\textbf{Theorem 1.2.} \emph{Let $F_t:\ M^n\longrightarrow
\mathbb{R}^{n+1}$ be a solution of the mean curvature flow of closed
hypersurfaces on a finite time interval $[0,T)$. If}\\
(1) \emph{there is a positive constant $C$ such that $h_{ij}\geq -C$
for $(x,t)\in M\times [0,T)$},\\
(2) \emph{$|| H||_{\alpha,M\times [0,T)}=\left(\int^T_0\int_M | H
|^\alpha d\mu dt\right)^{\frac{1}{\alpha}}<+\infty$ for some
$\alpha\geq n+2$,\\
then this flow can be extended over time $T$.}\\

When the initial hypersurface is mean convex, we have following\\\\
\textbf{Theorem 1.3.} \emph{Let $F_t:\ M^n\longrightarrow
\mathbb{R}^{n+1}$ be a solution of the mean curvature flow of closed
hypersurfaces on a finite time interval $[0,T)$. If}\\
(1) \emph{$H>0$ at $t=0$},\\
(2) \emph{$|| H||_{\alpha,M\times [0,T)}=\left(\int^T_0\int_M | H
|^\alpha d\mu dt\right)^{\frac{1}{\alpha}}<+\infty$ for some
$\alpha\geq n+2$,\\
then this flow can be extended over time $T$.}\\

The following example shows that the condition $\alpha\geq n+2$ in
Theorem 1.1, 1.2 and 1.3 is optimal.\\\\
\textbf{Example.} Set
$\mathbb{S}^n=\{x\in\mathbb{R}^{n+1}:\sum^{n}_{i=1}x_i^2=1\}$. Let
$F$ be the standard isometric embedding of $\mathbb{S}^n$ into
$\mathbb{R}^{n+1}$. It is clear that $F(t)=\sqrt{1-2nt}F$ is the
solution to the mean curvature flow, where $T=\frac{1}{2n}$ is the
maximal existence time. By a simple computation, we have
$g_{ij}(t)=(1-2nt)g_{ij}$, $H(t)=\frac{n}{\sqrt{1-2nt}}$ and
$h_{ij}(t)\geq0$. Hence
\begin{eqnarray*}
|| H||_{\alpha,M\times [0,T)}&=&\left(\int^T_0\int_M
| H |^\alpha d\mu dt\right)^{\frac{1}{\alpha}}\\
&=&C_1\left(\int^T_0
(T-t)^{\frac{n-\alpha}{2}}dt\right)^{\frac{1}{\alpha}},
\end{eqnarray*}
where $C_1$ is a positive constant. It follows that

\[|| H||_{\alpha,M\times
[0,T)}\left\{ \begin{array}{ll}
=\infty,&\ for\ \ \alpha \geq n+2,\\
<\infty,&\ for\ \ \alpha < n+2.
\end{array}
\right.\] This implies that the condition $\alpha\geq n+2$ in
Theorem 1.2 and Theorem 1.3 is optimal.

Since $F(t)$ is an umbilical hypersurface in $\mathbb{R}^{n+1}$ for
$t\in [0,T)$, $|| A||_{\alpha,M\times [0,T)}=\frac{1}{\sqrt{n}}||
H||_{\alpha,M\times [0,T)}$. Therefore,

\[|| A||_{\alpha,M\times
[0,T)}\left\{ \begin{array}{ll}
=\infty,&\ for\ \ \alpha \geq n+2,\\
<\infty,&\ for\ \ \alpha < n+2.
\end{array}
\right.\] So the condition $\alpha\geq n+2$ in Theorem 1.1 is also
optimal.

\section{Mean curvature flow with finite $L^{\alpha}$ norm of $A$}
In this section, we extend the mean curvature flow with finite
$L^{\alpha}$ norm of
the second fundamental form, and give the proof of Theorem 1.1.\\
\\
\emph{Proof of Theorem 1.1.}\ \ By H\"{o}lder's inequality, $||
A||_{\alpha,M\times [0,T)}<+\infty$ implies $|| A||_{n+2,M\times
[0,T)}<+\infty$ if $\alpha>n+2$. Thus it is sufficient for us to
prove the theorem in the case where $\alpha=n+2$.

We argue by contradiction.

Suppose that $T(<\infty)$ is the maximal existence time. Firstly we
choose a sequence of time $t^{(i)}$ such that $\lim_{i\rightarrow
\infty}t^{(i)}=T$. Then we take a sequence of points $x^{(i)}\in M$
satisfying
$$|A|^2(x^{(i)},t^{(i)})=\max_{(x,t)\in M\times [0,t^{(i)})}|
A|^2(x,t),\eqno(1)$$ where $\lim_{i\rightarrow \infty}|
A|^2(x^{(i)},t^{(i)})=+\infty$.

Putting $Q^{(i)}=| A|^2(x^{(i)},t^{(i)})$, we consider the rescaling
mean curvature flow:
$$F^{(i)}(t)=\left(Q^{(i)}\right)^{\frac{1}{2}}F\left(\frac{t}{Q^{(i)}}+t^{(i)}\right).\eqno(2)$$
Then the induced metric on $M$ by the immersion $F^{(i)}(t)$ is
$g^{(i)}(t)=Q^{(i)}g\left(\frac{t}{Q^{(i)}}+t^{(i)}\right)$, $t\in
(-Q^{(i)}t^{(i)},0]$. For $(M,g^{(i)}(t))$, the second fundamental
form $| A^{(i)}|(x,t)\leq 1$, for any $i$.

From \cite{2}, there exists a subsequence of
$(M,g^{(i)}(t),x^{(i)})$ that converges to a Riemannian manifold
$(\overline{M},\overline{g}(t),\overline{x})$, $t\in (-\infty, 0]$,
and the corresponding subsequence of immersions  $F^{(i)}(t)$
converges to an immersion $\overline{F}(t): \overline{M}\rightarrow
\mathbb{R}^{n+1}$. Since $\int^{T_2}_{T_1}\int_M|A|^{n+2}_{g(t)}d\mu
dt$ is invariant under the rescaling
$Q^{\frac{1}{2}}F(x,\frac{t}{Q})$ , we calculate that
\begin{eqnarray*}
\int^0_{-1}\int_{B_{\overline{g}(0)}(\overline{x},1)}|
\overline{A}|^{n+2}d\overline{\mu}dt&\leq &\lim_{i\rightarrow
\infty}\int^0_{-1}\int_{B_{g^{(i)}(0)}(x^{(i)},1)}|
A|^{n+2}_{g^{(i)}(t)}d\mu_{g^{(i)}(t)}dt\\
&=&\lim_{i\rightarrow
\infty}\int^{t^{(i)}}_{t^{(i)}-(Q^{(i)})^{-1}}\int_{B_{g(t^{(i)})}(x^{(i)},(Q^{(i)})^{-\frac{1}{2}})}|
A|^{n+2}d\mu dt\\
&\leq &\lim_{i\rightarrow
\infty}\int^{t^{(i)}}_{t^{(i)}-(Q^{(i)})^{-1}}\int_M| A|^{n+2}d\mu
dt\end{eqnarray*}
$$=\ \ 0.\ \ \ \ \ \ \ \ \ \ \ \ \ \ \ \ \ \ \ \ \ \ \ \ \ \eqno(3)$$
The equality in $(3)$ holds because $\int^T_0\int_M | A |^{n+2} d\mu
dt<+\infty$ and $\lim_{i\rightarrow \infty}(Q^{(i)})^{-1}=0$. Since
$(\overline{M},\overline{g}(t))$ is a smooth Riemannian manifold for
each $t\leq 0$, the equality in $(3)$ implies that
$|\overline{A}|\equiv 0$ on
$B_{\overline{g}(0)}(\overline{x},1)\times [-1,0]$. In particular,
$|\overline{A}|(\overline{x},0)=0$. However, the points selecting
process implies that
$$|\overline{A}|(\overline{x},0)=\lim_{i\rightarrow \infty}|
A|_{g^{(i)}}(x^{(i)},0)=1,\eqno(4)$$ which is a contradiction. This completes the proof of Theorem 1.1.\\

By a similar argument, we prove that if $\lim_{t\rightarrow
T}\left(\int_M | A |^\alpha d\mu\right
)^{\frac{1}{\alpha}}<+\infty$,
the mean curvature flow can be extended over time $T$.\\
\\
\textbf{Theorem 2.1.} \emph{Let $F_t:\ M^n\longrightarrow
\mathbb{R}^{n+1}$ be a solution of the mean curvature flow of closed
hypersurfaces on a finite time interval $[0,T)$. If
$$\lim_{t\rightarrow T}\left(\int_M | A |^\alpha d\mu\right )^{\frac{1}{\alpha}}<+\infty$$
for some $\alpha\geq n$, then this flow can be extended over time
$T$.}\\
\emph{Proof.}\ \ It is sufficient for us to prove the theorem in the
case where $\alpha=n$, and we argue by contradiction again.

Suppose that $T(<\infty)$ is the maximal existence time. Let
$(x^{(i)},t^{(i)})$, $Q^{(i)}$, $F^{(i)}(t)$, $g^{(i)}(t)$,
$A^{(i)}$ and $(\overline{M},\overline{g}(t),\overline{x})$ be the
same as in the proof of Theorem 1.1. Since $\int_M|A|^{n}d\mu$ is
invariant under the rescaling $QF(x,t)$, we have
\begin{eqnarray*}
\int_{B_{\overline{g}(0)}(\overline{x},1)}|
\overline{A}|_{\overline{g}(0)}^{n}d\overline{\mu}_{\overline{g}(0)}dt&\leq
&\lim_{i\rightarrow \infty}\int_{B_{g^{(i)}(0)}(x^{(i)},1)}|
A|^{n}_{g^{(i)}(t)}d\mu_{g^{(i)}(t)}\\
&=&\lim_{i\rightarrow
\infty}\int_{B_{g(t^{(i)})}(x^{(i)},(Q^{(i)})^{-\frac{1}{2}})}|
A|_{g(t^{(i)})}^{n}d\mu_{g(t^{(i)})}\end{eqnarray*}
$$=\ \ 0.\ \ \ \ \ \ \ \ \ \ \ \ \ \ \ \ \eqno(5)$$
The equality in $(5)$ holds because $\lim_{t\rightarrow
T}\left(\int_M | A |^\alpha d\mu\right )^{\frac{1}{\alpha}}<+\infty$
and
$B_{g(t^{(i)})}(x^{(i)},(Q^{(i)})^{-\frac{1}{2}})\rightarrow\emptyset$
as $i\rightarrow \infty$. The equality $(5)$ implies that
$|\overline{A}|\equiv 0$ on the ball
$B_{\overline{g}(0)}(\overline{x},1)$. In particular,
$|\overline{A}|(\overline{x},0)=0$. On the other hand, the points
selecting process implies that
$$|\overline{A}|(\overline{x},0)=\lim_{i\rightarrow \infty}|
A|_{g^{(i)}}(x^{(i)},0)=1.\eqno(6)$$The contradiction completes the
proof. It is easy to check that the  the condition $\alpha\geq n$ is
optimal.

\section{Mean curvature flow with finite total mean curvature}

In this section we prove Theorem 1.2. We first recall
some evolution equations (see [3], [14]).\\\\
\textbf{Lemma 3.1.} \emph{Let $g=\{g_{ij}\}$ and $A=\{h_{ij}\}$ be
the metric and the second fundamental form on $M$, and denote by
$H=g^{ij}h_{ij}$, $| A|^2=h^{ij}h_{ij}$ the mean curvature and the
squared norm of the second fundamental form respectively, then}
\begin{eqnarray*}
&& \frac{\partial}{\partial t}g_{ij}=-2Hh_{ij},\\
&& \frac{\partial \overrightarrow{\nu}}{\partial t}=\nabla^i
H\frac{\partial F}{\partial x^i},\\
&& \frac{\partial}{\partial t}h_{ij}=\triangle
h_{ij}-2Hh_{il}g^{lm}h_{mj}+| A|^2h_{ij},\\
&&\frac{\partial}{\partial t}H=\triangle H+| A|^2H,\\
&& \frac{\partial}{\partial t}| A|^2=\triangle| A|^2-2| \nabla
A|^2+2| A|^4.
\end{eqnarray*}

The following Sobolev inequality can be found in [8] and [12].\\\\
\textbf{Lemma 3.2.} \emph{Let $M$ be an $n-$dimensional $(n\geq 3)$
closed submanifold of a Riemannian manifold $N^{n+p}$ with
codimension $p\geq 1$. Suppose that the sectional curvature of
$N^{n+p}$ is non-positive. Then for any $s\in (0,+\infty)$ and $f\in
C^1(M)$ such that $f\geq 0$,
$$\int_M |\nabla f|^2\geq \frac{(n-2)^2}{4(n-1)(1+s)}\left[\frac{1}{C^2(n)}\left(\int_M f^{\frac{2n}{n-2}}\right)^{\frac{n-2}{n}}
-H^2_0\left(1+\frac{1}{s}\right)\int_M f^2\right],$$ where
$H_0=\max_{x\in M}|H|$,
$C(n)=\frac{2^n(1+n)^{(1+\frac{1}{n})}}{(n-1)\sigma_n}$, and
$\sigma_n$ is the volume of the unit ball in $\mathbb{R}^{n+1}$.}\\

The following estimate is very useful in the proof of Theorem 1.2.\\\\
\textbf{Theorem 3.3.} \emph{Suppose that $F_t: M\longrightarrow
\mathbb{R}^{n+1}$ is a mean curvature flow solution for $t\in
[0,T_0]$, and the second fundamental form is uniformly bounded on
time interval $[0,T_0]$. Then
$$\max_{(x,t)\in M\times [\frac{T_0}{2},T_0]} H^2(x,t)\leq C_2\left(\int^{T_0}_0\int_{M_t}
 | H|^{n+2}d\mu dt\right)^{\frac{2}{n+2}},$$
where $C_2$ is a constant depending on $n$, $T_0$ and
$\sup_{(x,t)\in M\times[0,T_0]}|A|$.}\\
\emph{Proof.} The evolution equation of $H^2$ is
$$\frac{\partial}{\partial t}H^2=\triangle H^2-2 |\nabla H|^2+2| A|^2H^2.\eqno(7)$$
Since $| A|$ is bounded, we obtain the following estimate from (7).
$$\frac{\partial}{\partial t}H^2\leq \triangle H^2+\beta H^2,\eqno(8)$$
where $\beta$ is a constant depending only on $\sup_{(x,t)\in
M\times[0,T_0]}|A|$.

Denoting $f=H^2$, from the inequality in $(8)$ we obtain that for
any $p\geq 2$,
\begin{eqnarray*}
\frac{\partial}{\partial t}\int_{M_t} f^{p}&=&\int_{M_t}pf^{p-1}\frac{\partial}{\partial t}f-\int_{M_t} f^{p+1} \\
&\leq&\int_{M_t} pf^{p-1}\left(\triangle f+\beta f\right)\\
&=&-\frac{4(p-1)}{p}\int_{M_t}|\nabla f^{\frac{p}{2}}|^2+\beta
p\int_{M_t} f^{p}. \end{eqnarray*} Thus
$$\frac{\partial}{\partial
t}\int_{M_t} f^{p}+\frac{4(p-1)}{p}\int_{M_t}|\nabla
f^{\frac{p}{2}}|^2\leq \beta p\int_{M_t} f^{p}.\eqno(9)$$

For any $0<\tau<\tau'<T_0$, define a function $\psi$ on $[0,T_0]$:\\
\[\psi(t)=\left\{ \begin{array}{ll}
0&\ \ \ \ \ \ 0\leq t\leq \tau,\\
\frac{t-\tau}{\tau'-\tau}&\ \ \ \ \ \ \tau\leq t\leq \tau',\\
1&\ \ \ \ \ \ \tau'\leq t\leq T_0.
\end{array}
\right.\] Then by (9) we have
\begin{eqnarray*}
\frac{\partial}{\partial t}\left(\psi\int_{M_t}
f^{p}\right)&=&\psi'\int_{M_t} f^{p}+\psi\frac{\partial}{\partial
t}\left(\int_{M_t} f^{p}\right)\ \ \ \ \ \ \ \ \ \ \ \ \ \ \ \ \ \ \
\ \ \ \ \ \ \ \ \ \ \ \ \ \ \
\end{eqnarray*}
$$\ \ \ \ \ \ \ \ \ \ \ \ \ \ \ \ \ \ \ \ \ \ \ \ \ \ \ \ \ \ \ \ \leq\ \ \psi'\int_{M_t}
f^{p}+\psi\left(-\frac{4(p-1)}{p}\int_{M_t}|\nabla
f^{\frac{p}{2}}|^2+\beta p\int_{M_t} f^{p}\right).\eqno(10)$$ For
any $t\in [\tau',T_0]$, integrating both sides of the inequality in
$(10)$ on $[\tau,t]$ we get
$$\int_{M_t} f^{p}+\frac{4(p-1)}{p}\int^t_{\tau'}\int_{M_t}|\nabla f^{\frac{p}{2}}|^2\leq
\left(\beta+\frac{1}{\tau'-\tau}\right)\int^{T_0}_{\tau}\int_{M_t}
f^{p}.\eqno(11)$$ For the integral $\int^{T_0}_{\tau'}\int_{M_t}
f^{p(1+\frac{2}{n})}$, by Schwarz inequality and Sobolev inequality
in Lemma 3.2, we have
\begin{eqnarray*}
\int^{T_0}_{\tau'}\int_{M_t} f^{p(1+\frac{2}{n})}&\leq
&\int^{T_0}_{\tau'}\left(\int_{M_t} f^{p}\right)^{\frac{2}{n}}
\left(\int_{M_t} f^{\frac{np}{n-2}}\right)^{\frac{n-2}{n}} \\
&\leq&\max_{t\in [\tau',T_0]}\left(\int_{M_t}
f^{p}\right)^{\frac{2}{n}}
\int^{T_0}_{\tau'}\left(\int_{M_t} f^{\frac{np}{n-2}}\right)^{\frac{n-2}{n}}\\
&\leq&\left(\beta+\frac{1}{\tau'-\tau}\right)^{\frac{2}{n}}\left(\int^{T_0}_{\tau}\int_{M_t}
f^{p}\right)^{\frac{2}{n}}\\
&&\times
\int^{T_0}_{\tau'}\left[\frac{4(n-1)C^2(n)(1+s)}{(n-2)^2}\int_{M_t}|\nabla
f^{\frac{p}{2}}|^2 +\frac{n}{2}\beta
C^2(n)\left(1+\frac{1}{s}\right)\int_{M_t}f^{p}\right].\end{eqnarray*}
For the third factor on the right hand side, we have from (11)
\begin{eqnarray*}
&&\int^{T_0}_{\tau'}\left[\frac{4(n-1)C^2(n)(1+s)}{(n-2)^2}\int_{M_t}|\nabla
f^{\frac{p}{2}}|^2+\frac{n}{2}\beta
C^2(n)\left(1+\frac{1}{s}\right)\int_{M_t}f^{p}\right]\\
&\leq&\frac{4(n-1)C^2(n)(1+s)}{(n-2)^2}\int^{T_0}_{\tau'}\int_{M_t}|\nabla
f^{\frac{p}{2}}|^2\\
&&+\frac{n}{2}\beta
C^2(n)\left(1+\frac{1}{s}\right)\int^{T_0}_{\tau'}\left[\left(\beta
p+\frac{1}{\tau'-\tau}\right)\int^{T_0}_{\tau}\int_{M_t}f^{p}\right]\\
&\leq&\frac{(n-1)C^2(n)p(1+s)}{(n-2)^2(p-1)}\left(\beta
p+\frac{1}{\tau'-\tau}\right)\int^{T_0}_{\tau}\int_{M_t}f^{p}\\
&&+\frac{n}{2}\beta C^2(n)T_0\left(1+\frac{1}{s}\right)\left(\beta
p+\frac{1}{\tau'-\tau}\right)\int^{T_0}_{\tau}\int_{M_t}f^{p}\\
&=&\left[\frac{(n-1)C^2(n)p(1+s)}{(n-2)^2(p-1)}+\frac{n}{2}\beta
C^2(n)T_0\left(1+\frac{1}{s}\right)\right]\left(\beta
p+\frac{1}{\tau'-\tau}\right)\int^{T_0}_{\tau}\int_{M_t}f^{p}.\end{eqnarray*}
Hence
\begin{eqnarray*}
&&\int^{T_0}_{\tau'}\int_{M_t}f^{p\left(1+\frac{2}{n}\right)}\\
&\leq&\left[\frac{(n-1)C^2(n)p(1+s)}{(n-2)^2(p-1)}+\frac{n}{2}\beta
C^2(n)T_0\left(1+\frac{1}{s}\right)\right]\left(\beta
p+\frac{1}{\tau'-\tau}\right)^{1+\frac{2}{n}}\left(\int^{T_0}_{\tau}\int_{M_t}f^{p}\right)^{1+\frac{2}{n}}.
\end{eqnarray*}
Put $L(p,t)=\int^{T_0}_{t}\int_{M_t}f^{p}$,
$s=\frac{\left[\frac{2}{n}(p-1)T_0\beta\right]^{\frac{1}{2}}(n-2)}{[n(p-1)]^{\frac{1}{2}}}$,
and $D=\frac{[(n-1)p]^{\frac{1}{2}}C(n)}{(n-2)(p-1)^{\frac{1}{2}}}$.
The above inequality can be rewritten as
$$L\left(p\left(1+\frac{2}{n}\right),\tau'\right)\leq D^2
\left(\beta
p+\frac{1}{\tau'-\tau}\right)^{1+\frac{2}{n}}L(p,\tau)^{1+\frac{2}{n}}.\eqno(12)$$
Now let $\mu=1+\frac{2}{n}$, $p_k=\frac{n+2}{2}\mu^k$ and
$\tau_k=\left(1-\frac{1}{\mu^{k+1}}\right)t$. Then from (12) we
obtain
$$L(p_{k+1},\tau_{k+1})^{\frac{1}{p_{k+1}}}\leq D^{\Sigma_{i=0}^{k}\frac{2}{p_{i+1}}}\left(\frac{(n+2)\beta}{2}+
\frac{n+2}{2t}\right)^{\Sigma_{i=0}^{k}\frac{1}{p_i}}\mu^{\Sigma_{i=0}^{k}\frac{i}{p_i}}L(p_0,\tau_0)^{\frac{2}{n+2}}.$$
As $k\rightarrow +\infty$, we conclude
$$f(x,t)\leq D^{\frac{2n}{n+2}}\left(1+\frac{2}{n}\right)^{\frac{n}{2}}\left(\frac{n+2}{2}\beta+\frac{n+2}{2t}\right)
\left(\int^{T_0}_0\int_{M_t}f^{\frac{{n+2}}{2}}\right)^{\frac{2}{n+2}}.\eqno(13)$$
Therefore, for any $(x,t)\in M\times[\frac{T_0}{2},T_0]$, we get
from (13)
$$H^2(x,t)\leq C_2\left(\int^{T_0}_0\int_{M_t}| H|^{n+2}\right)^{\frac{2}{n+2}},$$
where $C_2$ is a constant depending on $n$, $T_0$ and
$\sup_{(x,t)\in M\times[0,T_0]}|A|$. Thus
$$\max_{(x,t)\in M\times[\frac{T_0}{2},T_0]}H^2(x,t)\leq C_2\left(\int^{T_0}_0\int_{M_t}| H|^{n+2}d\mu dt\right)^{\frac{2}{n+2}},$$
which is desired.\\

We are now in a position to prove Theorem 1.2.\\\\
\emph{Proof of Theorem 1.2.}\ \ We only need to prove the theorem
for $\alpha=n+2$ since by H\"{o}lder's inequality, $||
H||_{\alpha,M\times [0,T)}<\infty$ implies $|| H||_{n+2,M\times
[0,T)}<\infty$ if $\alpha>n+2$. We still argue by contradiction.

Suppose that the solution to the mean curvature flow can't be
extended over $T$, then $| A|$ becomes unbounded as $t\rightarrow
T$. Since $h_{ij}\geq -C$, we get $\sum_{i,j}(h_{ij}+C)^2\leq
C_3[tr(h_{ij}+C)]^2$, where $C_3$ is a constant depending only on
$n$. On one hand, $| A|^2$ is unbounded implies that
$\sum_{i,j}(h_{ij}+C)^2$ is unbounded. On the other hand,
$$[tr(h_{ij}+C)]^2=(H+nC)^2=H^2+2nCH+n^2C^2.\eqno(14)$$
Thus $H^2$ is unbounded. Namely, $$\sup_{(x,t)\in M\times
[0,T)}H^2(x,t)=\infty.$$ Choose an increasing time sequence
$t^{(i)}$, $i=1,2,\cdots$, such that $\lim_{i\rightarrow
\infty}t^{(i)}=T$. We take a sequence of points $x^{(i)}\in M$
satisfying
$$H^2(x^{(i)},t^{(i)})=\max_{(x,t)\in M\times [0,t^{(i)})} H^2(x,t).$$
Then $\lim_{i\rightarrow \infty} H^2(x^{(i)},t^{(i)})=\infty$.

Putting $Q^{(i)}=H^2(x^{(i)},t^{(i)})$, we have
$\lim_{i\rightarrow\infty}Q^{(i)}=\infty$. This together with
$\lim_{i\rightarrow \infty}t^{(i)}=T>0$ implies that there exists a
positive integer $i_0$ such that $Q^{(i)}t^{(i)}\geq1$ for $i\geq
i_0$. For $i\geq i_0$ and $t\in [0,1]$, we define
$F^{(i)}(t)=\left(Q^{(i)}\right)^{\frac{1}{2}}F\left(\frac{t-1}{Q^{(i)}}+t^{(i)}\right)$.
Then the induced metric on $M$ induced by $F^{(i)}(t)$ is
$g^{(i)}(t)=Q^{(i)}g\left(\frac{t-1}{Q^{(i)}}+t^{(i)}\right)$, and
$F^{(i)}(t):M\rightarrow \mathbb{R}^{n+1}$ is still a solution to
the mean curvature flow on $t\in[0,1]$. Since $F_t$ satisfies
$h_{ij}\geq-C$ for $(x,t)\in M\times[0,T)$, we have
$$H^2_{(i)}(x,t)\leq 1 \ \ on\ \ M\times [0,1] ,$$
$$h^{(i)}_{jk}\geq -\frac{C}{\sqrt{Q^{(i)}}}\ \ on\ \ M\times [0,1],\eqno(15)$$
where $H_{(i)}$ and $A^{(i)}=h^{(i)}_{jk}$ are mean curvature and
the second fundamental form of $F^{(i)}(t)$ respectively. The
inequality in $(15)$ gives that
$h^{(i)}_{jk}+\frac{C}{\sqrt{Q^{(i)}}}\geq 0$. Hence
$$h^{(i)}_{jk}+\frac{C}{\sqrt{Q^{(i)}}}\leq tr\left(h^{(i)}_{jk}+\frac{C}{\sqrt{Q^{(i)}}}\right)
\leq H_{(i)}+\frac{nC}{\sqrt{Q^{(i)}}},\eqno(16)$$ which implies
that $h^{(i)}_{jk}\leq H_{(i)}+\frac{(n-1)C}{\sqrt{Q^{(i)}}}$. Also,
since $Q^{(i)}\rightarrow +\infty$ as $i\rightarrow \infty$, we know
that $h^{(i)}_{jk}\leq C_4$, where $C_4$ is a constant independent
of $i$.

Set
$(M^{(i)},g^{(i)}(t),x^{(i)})=\left(M,Q^{(i)}g\left(\frac{t-1}{Q^{(i)}}+t^{(i)}\right),x^{(i)}\right)$,
$t\in [0,1]$. From \cite{2} we know that there is a subsequence of
$(M^{(i)},g^{(i)}(t),x^{(i)})$ converges to a Riemannian manifold
$(\widetilde{M},\widetilde{g}(t),\widetilde{x})$, and the
corresponding subsequence of immersions $F^{(i)}(t)$ converges to an
immersion $\widetilde{F}(t):\widetilde{M}\rightarrow
\mathbb{R}^{n+1}$.

Since $F^{(i)}(t)$ satisfies $H^2_{(i)}\leq 1$ on $M\times[0,1]$ for
any $i\geq i_0$, we know that $A^{(i)}$ is bounded by a constant
independent of $i$, for $t\in [0,1]$. It follows from Theorem 3.3
that
$$\max_{(x,t)\in M^{(i)}\times[\frac{1}{2},1]}H^2_{(i)}(x,t)\leq C_5\left(\int^{1}_0\int_M| H|^{n+2}_{(i)}d\mu_{g^{(i)}(t)}
dt\right)^{\frac{2}{n+2}},$$ where $C_5$ is a constant independent
of $i$. Since $\int^{T_2}_{T_1}\int_M|H|^{n+2}_{g(t)}d\mu dt$ is
invariant under the rescaling $Q^{\frac{1}{2}}F(x,\frac{t}{Q})$,
using similar calculation as in (3) we have
\begin{eqnarray*}
\max_{(x,t)\in
\widetilde{M}\times[\frac{1}{2},1]}\widetilde{H}^2(x,t)&\leq
&\lim_{i\rightarrow \infty}C_5\left(\int^{1}_0\int_M|
H|^{n+2}_{(i)}d\mu_{g^{(i)}(t)}
dt\right)^{\frac{2}{n+2}}\\
&\leq&\lim_{i\rightarrow
\infty}C_5\left(\int^{t^{(i)}+(Q^{(i)})^{-1}}_{t^{(i)}}\int_M|
H|^{n+2}_{(i)}d\mu dt\right)^{\frac{2}{n+2}}\end{eqnarray*}
$$=\ \ 0.\ \ \ \ \ \ \ \ \ \ \ \ \ \ \ \ \ \ \ \ \ \ \ \ \ \ \eqno(17)$$
The equality in $(17)$ holds because $\int^T_0\int_M  H^{n+2} d\mu
dt<+\infty$ and $\lim_{i\rightarrow \infty}(Q^{(i)})^{-1}=0$.\\
However, according to the choice of the points, we have
$$\widetilde{H}^2(\widetilde{x},1)=\lim_{i\rightarrow \infty}H^2_{(i)}(x^{(i)},1)=1.$$
This is a contradiction. We complete the proof of Theorem 1.2.\\

With a similar method, we can prove Theorem 1.3.\\\\
\emph{Proof of Theorem 1.3.} Since $H>0$ at $t=0$, there exists a
positive constant $C_6$ such that $|A|^2\leq C_6 H^2$. The evolution
of $H$ in Lemma 3.1 implies that $H>0$ is preserved along the mean
curvature flow. By [7] we have the following evolution equation of
$\frac{|A|^2}{H^2}$:
$$\frac{\partial}{\partial t}\left(\frac{|A|^2}{H^2}\right)=\triangle\left(\frac{|A|^2}{H^2}\right)+\frac{2}{H}
\left\langle\nabla
H,\nabla\left(\frac{|A|^2}{H^2}\right)\right\rangle-\frac{2}{H^4}|H\nabla_ih_{jk}-\nabla_iH\cdot
h_{jk}|^2.\eqno(18)$$ From the maximum principle, we obtain that
$|A|^2\leq C_6 H^2$ is preserved along the mean curvature flow.

It is sufficient to prove the theorem for $\alpha=n+2$. We still
argue by contradiction. Suppose that the solution to the mean
curvature flow can't be extended over time $T$, then  $|A|^2$ is
unbounded as $t\rightarrow T$. This implies that $H^2$ is also
unbounded since $|A|^2\leq C_6 H^2$. Let $(x^{(i)},t^{(i)})$,
$Q^{(i)}$, $F^{(i)}(t)$, $g^{(i)}(t)$ and
$(\widetilde{M},\tilde{g}(t),\tilde{x})$ be the same as in the proof
of Theorem 1.2. Let $A^{(i)}$ and $H_{(i)}$ be the second
fundamental form and mean curvature of the immersion $F^{(i)}(t)$
respectively. Then we have $|A^{(i)}|^2\leq C_6|H_{(i)}|^2$ for
$(x,t)\in M\times [0,1]$, which implies that $A^{(i)}$ is bounded by
a constant independent of $i$, for $t\in [0,1]$. It follows from
Theorem 3.3 that
$$\max_{(x,t)\in M^{(i)}\times[\frac{1}{2},1]}H^2_{(i)}(x,t)\leq C_7\left(\int^{1}_0\int_M| H|^{n+2}_{(i)}d\mu_{g^{(i)}(t)}
dt\right)^{\frac{2}{n+2}},$$ where $C_7$ is a constant independent
of $i$. By an argument similar to the proof of Theorem 1.2, we can
get a contradiction which completes the proof of Theorem 1.3.\\

Finally we would like to propose the following\\
\textbf{Open question}.\ \ Can one generalize Theorems 1.1, 1.2 and
1.3 to the case where the ambient space is a general Riemannian
manifold?



Center of Mathematical Sciences\

Zhejiang University

Hangzhou 310027

China\\

E-mail address: xuhw@cms.zju.edu.cn; yf@cms.zju.edu.cn;
superzet@163.com

\end{document}